\newcommand{\di}{{\: \rm d}}   

\newcounter{secti}
\newcounter{item}[secti]
\renewcommand{\theitem}{\thesecti.\arabic{item}}

\newcommand{\sect}[1]{\medskip \begin{center} \refstepcounter{secti}
{\sc \thesecti.\  #1 }\end{center} \medskip}

\newcommand{\qitem}[2]{\refstepcounter{item}
  {\bf \theitem\  #1.}\  {#2 }}

\newcommand{\eitem}{\medskip}
\newcommand{\cF}{{\mathcal F}}

\newcommand{\Leb}{\, {\rm Leb\, }}
\newcommand{\lee}[1]{#1_{\ast}}
\newcommand{\ree}[1]{#1^{\ast}}

\newcommand{\R}{\mathbb{R}}

\newcommand{\be}{\begin{equation} }

\newcommand{\ee}{\end{equation}}
\newcommand{\beq}{\begin{eqnarray} }
\newcommand{\eeq}{\end{eqnarray}}
\newcommand{\beqn}{\begin{eqnarray*} }
\newcommand{\eeqn}{\end{eqnarray*}}
\newcommand{\tends}{\rightarrow}
\newcommand{\kraj}{$\quad\Box$}
\newcommand{\sgn}{{\rm sgn\; }}

\newcommand{\E}{{\rm E\, }} 
\newcommand{\Med}{{\rm Med\; }} 

\documentclass{amsart}

\begin{document}

\begin{center}

{\LARGE\bf  Young's functional with Lebesgue-Stieltjes integrals}

\vspace{3em}

{\sc  Milan Merkle, Dan Marinescu, Monica M. R.  Merkle, Mihai Monea and Marian Stroe}

\vspace{3em}

\parbox{25cc}{
{\bf Abstract. }{\small For non-decreasing real  functions $f$ and $g$, we consider the functional
$ T(f,g ; I,J)=\int_{I} f(x)\di g(x) + \int_J g(x)\di f(x)$, where $I$ and $J$ are intervals with
$J\subseteq I$.  In particular case with $I=[a,t]$, $J=[a,s]$, $s\leq t$ and $g(x)=x$, this reduces to
the expression in classical Young's inequality. We  survey some properties of Lebesgue-Stieltjes interals
and present a new simple proof for change of variables. Further, we formulate a version of
Young's inequality with respect
to arbitrary positive finite measure on real line including a purely discrete case, and discuss
an application related to medians of probability distributions and a summation formula that involves values of a function
and its inverse at integers.

}
}

\end{center}

\medskip


\medskip


\vspace{3em}

\sect{Introduction and first generalization}\label{intro}

Let $f$ be a continuous and increasing function defined on an interval $I\subset \R$  and let $a,s,t$ be arbitrary points in $I$.
A version of classical Young's inequality \cite{young12} states that
\be
\label{yng1}
Y(f;a,s,t):= \int_a^t f(x)\di x +\int_{f(a)}^{f(s)} f^{-1}(x)\di x \geq tf(s)- af(a).
\ee

There has been a considerable amount of literature related to this classical inequality.
It seems that the first strict analytic proofs were given in \cite{diazmet70} and in \cite{bullen71}.
A relation between Young's functional and integration by parts in Riemann-Stieltjes integral was used in
 in \cite{boasmar74a} and \cite{boasmar74b}
for proving various versions of Young's inequality, and also in \cite{boasmar73} to derive several interesting applications.
 The reverse, in the sense
that if (\ref{yng1}) holds for some functional $Y(f;a,s,t)$ and for all continuous increasing $f$, then $Y$ is of the form as in (\ref{yng1}),
has also been a topic of research. Although in this article we do not address the question of reverse,
some interesting contributions
like \cite{bullen71}, \cite{hsu72} or \cite{lack74} are worth mentioning. Finally, let us recall that the topic is   closely related  to convexity,
complementary Young's functions and  Orlicz spaces.

 The first recorded result on upper bounds for Young's functional $Y$ was found by
 Merkle \cite{mer74} in 1974;  in the setup of (\ref{yng1}) it reads

\be
\label{merub}
Y(f;a,s,t) \leq  f(a)(t-a)+ a(f(s)-f(a))
 +\max ((t-a)(f(t)-f(a)), (s-a)(f(s)-f(a))),
\ee
where $s>a,t>a$. In 2008, Minguzzi \cite{ming08} discovered a closer  upper bound  also with $s>a,t>a$:
\be
\label{ming1}
Y(f;a,s,t) \leq tf(t) +sf(s)-af(a)-sf(t)
\ee
and observed that, given $f$ and an $a\in I$,  (\ref{merub}) holds for all $s,t>a$
if  and only if (\ref{ming1}) does, that is, the statements (\ref{merub}) and (\ref{ming1}) are equivalent. Incidentally, the same result was
stated  by Witkowski in \cite[Second proof of Theorem 1]{witk06} as well as
(together with the observation about equivalence) by Cerone
 \cite{cerone09} in 2009.

Extending the idea of \cite{boasmar74b}, one can observe a more general functional:

\beq
 \label{yngs}
 \lefteqn{S(f,g;a,s,t) :=  \int_a^tf(x)\di g(x) +\int_a^s g(x)\di f(x)}\\
 & & =\int_a^tf(x)\di g(x) +\int_{f(a)}^{f(s)} g(f^{-1}(y))\di y,\nonumber
  \eeq
  where $g$ is an arbitrary continuous and increasing function and integrals are taken in the Riemann-Stieltjes sense. For
  $g(x)=x$ this reduces to the classical case (\ref{yng1}). The next theorem presents our
  first generalization of (\ref{yng1}) and (\ref{ming1}).

\qitem{Theorem}{\label{ygen1}\it Let $f$  and $g$ be  continuous increasing functions  defined on an interval $I\subset \R$  and
let $a,s,t$ be arbitrary points in $I$. Then

\beq
 \label{yngss}
 g(t)f(s)-g(a)f(a) &\leq & S(f,g;a,s,t)\\
                  & \leq & g(t)f(t)+g(s)f(s)-g(a)f(a)-g(s)f(t),\nonumber
  \eeq
  with equality in both sides if and only if $s=t$.}

{\bf Proof. } Using integration by parts, we have that $S(f,g;a,s,s)= f(s)g(s)-f(a)g(a)$ and so, for $s\leq t$,
 \beqn
 \lefteqn{S(f,g;a,s,t) = S(f,g; a, s,s) + \int_s^t f(x)\di g(x)} \\
            &\geq & f(s)g(s)-f(a)g(a) + f(s)(g(t)-g(s))= g(t)f(s)-g(a)f(a),
 \eeqn
 which is the first inequality in (\ref{yngss}). The equality is possible if an only if
\[ \int_s^t (f(x)-f(s))\di g(x) =0, \]
  that is, if and only if $s=t$. For $s>t$, the proof is based on the relation
 \[ S(f,g;a,s,t) = S(f,g; a, t,t) + \int_t^s g(x)\di f(x) \]
 and then one proceeds as above.

 For the second inequality in (\ref{yngss}), we start with the observation that
 \[ S(f,g; a,s,t) + S(g,f; a,s,t) = f(t)g(t)-f(a)g(a)+f(s)g(s)-f(a)g(a) \]
 and so,
 \[ S(f,g; a,s,t)  = f(t)g(t)+f(s)g(s)-2f(a)g(a) - S(g,f; a,s,t).\]
 Applying here the first inequality in (\ref{yngss}) with $S(g,f;a,s,t)$ we get the second one. \kraj

 Let us remark that left and right inequality in (\ref{yngss}) follow from each other, so the equivalence between (\ref{yng1}) and (\ref{ming1})
 is preserved in the above generalization.

In section 3 we offer an ultimate generalization, for the case when $f$ and $g$ are non-increasing and not necessarily continuous.
To achieve this goal, we need the material of Section 2. Examples of applications are postponed to Section 4.

\eitem

\sect{Lebesgue-Stieltjes integral:  change of variables and integration by parts\label{lsip} }

In the proof of Theorem \ref{ygen1} we used integration by parts formula, and the representation in the second line of (\ref{yngs})
is due to change of variables $y=f(x)$. These are  extensions of formulas in elementary calculus to Riemann-Stieltjes case, and proofs
can be found
in the classical text \cite{rienagy}. In this section we deal with  Lebesgue-Stieltjes integrals of the form $\int f\di g$, where $f$ and $g$
are non-decreasing functions that may have jumps. The change of variables in this case is based on the notion of generalized inverse (defined
bellow, see (\ref{geninv})) and we were not able to locate a strict proof in the existing literature,
except a brief note in \cite[p.124]{rienagy}. Here we provide a simple proof (Lemma \ref{lscv}) for our case  with
$f$ and $g$ being non-decreasing; it can be extended to a general $f$, but this will be addressed elsewhere. The second part of this
section is devoted to integration by parts in Lebesgue-Stieltjes integrals.

For a non-decreasing function $g: \R \mapsto \R$ one can define a positive measure $\mu_g$ on the family of intervals $[a,b]$ by
\be
\label{mug}
\mu_g ([a,b]):= g(b_{+})-g(a_{-}),\qquad (a\leq b)
\ee
and extend it to the Borel sigma field on real line. The Lebesgue-Stieltjes (L-S) integral of a given measurable
function $f$ with respect to $g$ is
defined as the Lebesgue integral of $f$ with respect to $\mu_g$
\[
 \int_{\R} f(x)\di g(x) := \int_{\R}  f(x)\di \mu_g (x)
 \]
If this integral is finite, we say that $f$ is $g$-integrable. Let us assume that $f$ is   $g$-integrable and non-negative. Then from the
theory of Lebesque integration it is enough to consider a  partition of interval $[0,n]$ on the $y$-axis with subintervals
 $J_{n,k}= \left[\frac{k-1}{2^n}, \frac{k}{2^n}\right), k=1,2,\ldots, n2^n$, and points $y_{n,k} = \frac{k-1}{2^n}$. Then
 \be
 \label{limsum}
  \int_{\R} f(x)\di g(x) = \lim_{n\tends +\infty}\sum_{k=1}^{n2^n} y_{n,k} \mu_g (f^{-1}(J_{n,k})).
  \ee
In a general case, we use the decomposition $f=f_{+}-f_{-}$, where $f_{+}(x)= \max\{f(x),0\}$ and  $f_{-}(x)=\max\{-f(x),0\}$.

The equality (\ref{mug}) can be  related to Lebesgue measure on real line, since
$\mu_g ([a,b]) = g(b_{+})-g(a_{-}) = {\rm Leb\; } ([g(a_{-}), g(b_{+})])$. More generally, for an interval $I$ define the interval
 $I_g$ to be $[g(a_{-}),g(b_{+})]$,
 $[g(a_{-}),g(b_{-})]$, $[g(a_{+}),g(b_{+})]$, $[g(a_{+}),g(b_{-})]$, if $I$ is $[a,b]$, $[a,b)$, $(a,b]$ or $(a,b)$, respectively.
 Sometimes it is convenient to think about empty set as a special interval; for $I=\emptyset$ we define $I_g=[a,a]$
 where $a$ is any real number, say zero; then $\mu_g(\emptyset)={\rm Leb\; }([a,a])=0$.
In what follows we will need a uniform notation for ends of the closed interval $I_g$, so let us define $\lee{g}(I)$  and $\ree{g}(I)$ to be
left and right endpoint of the closed interval $I_g$. Note that for a non-decreasing $g$, $\lee{g}(I)\leq g(x)\leq \ree{g}(I)$ for any
$x\in I$. Then on an arbitrary interval $I$ we may define
\[ \int_I f(x)\di g(x) = \int_{\R} f(x)\di g_I(x),\]
where $g_I(x)$ is defined to be $g(x)$ if $\lee{g}(I)\leq g(x)\leq  \ree{g}(I)$, and $\lee{g}(I)$ or $\ree{g}(I)$ if $g(x)< \lee{g}(I)$ or
$g(x)> \ree{g}(I)$, respectively. If a non-decreasing $f$ is defined only on $I$, we may extend it to $\R$ in a way analogous to previous case,
and so integrals of the form $\int_I f(x)\di g(x)$ where $f$ and $g$ are non-decreasing functions can be always expressed via integrals
of the form $\int_{R}f(x)\di g(x)$ with some other non-decreasing functions $f$ and $g$ defined on $\R$. Therefore, we may always assume that
the functions $f$ and $g$ are defined  on $\R$.

Now, the general formula analogous to (\ref{mug}) reads
 \[
 \mu_g (I)= {\rm Leb\; }(I_g),\quad {\rm i.e.} \quad \int_I \di g(x) = \int_{I_g} \di y,
  \]
  which is in fact a change of variables $y=g(x)$.

Let  $g$ be a non-decreasing (not necessarily continuous) function on $I$. Following \cite{cugross71} and \cite{boasmar74b},
we say that $x$  a generalized inverse  of $g$ at $y$, $x=g^{-1}(y)$, if
\be
\label{geninv}
 \sup\{ t\; |\; g(t)<y\} \leq x \leq \inf\{ t\; |\; g(t)> y\} .
 \ee

The generalized inverse $g^{-1}(y)$ is not unique at a point $y$ if and only if $g(I)=\{y\}$, where $I$ is an interval
with non-empty interior; in that case, any rule that will assign an $x\in \bar{I}$ yields a version  of  $g^{-1}(y)$. From
(\ref{geninv}) the following properties follow easily:
\be
\label{propi1}
g(x) <y \implies  x\leq g^{-1} (y) , \quad g(x)>y \implies x\geq g^{-1}(y);
\ee
\be
\label{propi2}
 g(x_{-})<y<g(x_{+}) \implies g^{-1}(y)=x, \ \ \ g^{-1}(y)=x\implies g(x_{-})\leq y\leq g(x_{+});
 \ee
\be
\label{propi3}
 g^{-1}(y)< x \implies y\leq g(x_{-}),\quad g^{-1}(y)>x\implies y\geq g(x_{+});
 \ee
 \be
 \label{propi4}
 g^{-1}(y) \leq x \implies y \leq g(x_{+}),\quad g^{-1}(y) \geq x \implies y \geq g(x_{-}).
 \ee

\medskip

The next lemma gives another way of expressing $\mu_g(I)$, for an interval $I$.

\qitem{Lemma}{ \label{lemleb}\it Let   $g$ be a non-decreasing functions on $\R$ and let $I$ be an interval. For  any version of
the generalized inverse $g^{-1}$ it holds that $\{ y\; |\; g^{-1}(y) \in I\} \subseteq I_g$ and
\be
\label{lebeq}
 \mu_g(I)= {\rm Leb\; } (I_g) = {\rm Leb\; } (\{ y\; |\; g^{-1}(y)\in I \}).
\ee
}

{\bf Proof.}  We have to prove the second equality in (\ref{lebeq}) only.
For a given $I$ and a fixed version of $g^{-1}$, let $J:=\{ y\; |\; g^{-1}(y)\in I \}$.

 For a singleton  $I=\{a\}$, we have that $I_g= [g(a_{-}),g(a_{+})]$ and by the second property in (\ref{propi2}) we have that
$J \subseteq I_g $. The first property in (\ref{propi2}) implies that $I_g\subseteq J\cup \{g(a_{-}),g(a_{+})\}$, hence
$\Leb (I_g)=\Leb (J)$.

If $I$ is an open interval $I=(a,b)$ where $a<b$,  then $I_g= [g(a_{+}),g(b_{-})]$, and the property (\ref{propi3}) shows that
$J\subseteq I_g$, whereas the reverse of (\ref{propi4}) implies that $I_g\subseteq J\cup \{ g(a_{+}), g(b_{-})\}$, and again
$\Leb (I_g)=\Leb (J)$.

For other kinds of intervals, the statement follows from additivity. \kraj

\eitem

Now we can give  the announced proof of change of variables formula in case when $f$ is also non-decreasing.

\qitem{Lemma}{\label{lscv}\it Let $f$ and $g$ be non-decreasing functions on an interval $I$
 and let $g^{-1}$ be any version of generalized inverse. Then
\be
\label{e1gv}
 \int_{I} f(x)\di g(x) = \int_{I_g} f(g^{-1}(y))\di y .
 \ee
 }

{\bf Proof. }  Since $f$ is non-decreasing, then
 $f^{-1}(J)$ is an interval (or empty set) for any interval $J$. Then we use Lemma \ref{lemleb}
  to show that  terms of summation on the right hand side of  (\ref{limsum}) are equal for both integrals
  in  (\ref{e1gv}); therefore, the integrals are equal.   \kraj

 \eitem

The following result on integration by parts for  in L-S integrals  is due to
E. Hewitt \cite{hew60} (see also \cite[Theorem 6.2.2.]{lesti}).

\medskip

\qitem{Lemma}{\it \label{lsle} Let $f$ and $g$ be real valued  functions of finite variation defined on $\R$
and let $I$ be an interval.   If $D(I)$ is (at most countable) set of points
 $d\in I$   where  both $f$ and $g$ are  discontinuous, then
 \be
 \label{mgen}
\int_{I}f(x)\di g(x) +\int_{I} g(x)\di f(x)
               =  \mu_{fg}(I) +\sum_{d\in D(I)} A(d),
 \ee

where
\be
\label{defofa}
 A(d)= f(d)(g(d_{+})-g(d_{-})) + g(d)(f(d_{+})-f(d_{-}))
        -f(d_{+})g(d_{+}) +f(d_{-})g(d_{-}).
     \ee
}

Since L-S integral in explicit notations depends on order of points and type of intervals, the explicit form  of (\ref{mgen}) is different
for different types of intervals. For example, if $I=[a,b]$, $a<b$, then

\be
\label{mgenci}
 \int_{[a,b]}f(x)\di g(x) +\int_{[a,b]} g(x)\di f(x)
               =  f(b_{+})g(b_{+})-f(a_{-})g(a_{-}) +\sum_{d\in D([a.b])} A(d)
 \ee

In a special case when one of functions $f,g$ is right continuous at $d$ and both are non-decreasing (or both non-increasing)
\be
\label{deforc}
A(d)= (f(d)-f(d_{-}))(g(d)-g(d_{-})) \geq 0.
\ee

Let us finally remark that, given a finite measure $\mu$, we may define a right continuous function $g_r(x)=\mu(-\infty,x]$
so that $\mu=\mu_{g_r}$ in the sense of L-S measure.
Hence for an arbitrary positive finite measure $\mu$ and a non-decreasing function $f$, we may take $g=g_r$
to make summands $A(d)$ non-negative.

\newpage

\sect{Young's functional-lower and upper bounds \label{ylowup}}

In this section we  assume that all functions are defined on $\R$, as explained  in the previous section.
We also allow a possibility
that intervals of integration can be infinite, as far as  the value of the integral is finite. Theorem \ref{teodi} below
offers one possible version
of Young's inequality with Lebesque-Stieltjes integrals. Here we will use the following notations,
\[ \bar{f}(I)= \inf\{ f(x)\; |\; x>I\},\quad \underline{f}(I)=  \inf\{ f(x)\; |\; x>I\},\]
\[ f_{\min}(I)=\inf\{ f(x)\; |\; x\in I\}, \quad f_{\max}(I) =\sup\{ f(x)\; |\; x\in I\} ,\]
where $I$ is an arbitrary nonempty interval. We note that for any non-decreasing function $f$ we have
\[ \underline{f}(I) \leq \lee{f}(I) \leq f_{\min}(I) \leq f_{\max}(I)\leq \ree{f}(I)\leq \bar{f}(I).\]
In this section we present lower and upper bounds for the functional
\beq
\label{oper}
 \lefteqn{T(f,g ; I,J):=\int_{I} f(x)\di g(x) + \int_J g(x)\di f(x) }\\
 &=& \int_{I} f(x)\di g(x) + \int_{J_f} g(f^{-1}(u))\di u \nonumber
 \eeq
 where $I$ and $J$ are arbitrary intervals with $J\subseteq I$. For $I=J$, the functional (\ref{oper})  can be expressed via
 integration by parts formula  (\ref{mgen}).  A particular case of (\ref{oper}) that corresponds to classical Young's inequality is
 with $I=[a,t]$, $J=[a,s]$, $a<s\leq t$ and with $g(x)=x$.

\qitem{Theorem}{\label{teodi}\em Let $f$ and $g$ be non-decreasing functions on $\R$ and let $I,J$ be intervals such that $J\subseteq I$. Then

\beq
\label{youngg1}
T(f,g;I,J) &\geq & T(f,g;J,J) + f_{\min}(I)(\lee{g}(J)-\lee{g}(I)) \nonumber \\
 & & + \bar{f}(J) (\ree{g}(I) -\ree{g}(J)),
\eeq
with equality if $I=J$.

In particular, with $a<s\leq t$, we have that

\beq
\label{mres}
\lefteqn{\int_{[a,t]} f(x)\di g(x) +\int_{[a,s]} g(x)\di f(x)}\\
& \geq & f(s_{+})g(t_{+}) -f(a_{-})g(a_{-}) +\sum_{d\in D([a,s])} A(d),\nonumber
\eeq
where $A(d)$ is defined as in (\ref{defofa}).
}

{\bf Proof.} Let $J\subseteq I$. The set $I\setminus J$ is of the form $A\cup B$, where $A$ and $B$ are disjoint intervals or empty sets with
$a<x<b$ for any $a\in A, x\in J, b\in B$. We start from
\be
\label{p3}
T(f,g; I,J) = T(f,g; J,J) +  \int_{A\cup B} f(u)\di g(u)
\ee
and use a  bound
\beq
\label{simplebo}
\int_{A\cup B} f(u)\di g(u) &\geq &  f_{\min}(I)\mu_g(A)+\bar{f}(J)\mu_g(B)\\
 &=& f_{\min}(I)(\lee{g}(J)-\lee{g}(I)) + \bar{f}(J) (\ree{g}(I) -\ree{g}(J)),\nonumber
\eeq
where the second line follows from noticing that $\ree{g}(A)=\lee{g}(J), \lee{g}(A)=\lee{g}(I)$ and analogously for $B$.
Then (\ref{youngg1}) follows from (\ref{p3}) and (\ref{simplebo}).

The special case (\ref{mres}) follows from (\ref{youngg1}) for $I=[a,t]$, $J=[b,s]$, by using the integration by parts formula of Lemma \ref{lsle}
for $T(f,g; J,J)$.

\eitem

\qitem{Remark}{\label{simplupb} A simple observation that (\ref{simplebo}) can be re-phrased in terms
of upper bounds, yields an immediate upper bound in Young's
inequality.  The bound obtained in such a way agrees with the bound
that is provided in the next theorem if $f$ and $g$ are continuous functions.}

\eitem

\qitem{Theorem}{\label{teoco}\em Under same conditions and notations as in Theorem \ref{teodi},
\beq
\label{youngc}
T(f,g; I,J) &\leq & T(f,g; I,I) - g_{\min}(I) (\lee{f}(J)-\lee{f}(I)) \nonumber\\
& & - \bar{g}(J) (\ree{f}(I)-\ree{f}(J)),
\eeq

with equality if $I=J$. In particular, for $a\leq s\leq t$,

\beq
\label{mres1}
\lefteqn{\rule{7em}{0em}\int_{[a,t]} f(x)\di g(x) +\int_{[a,s]} g(x)\di f(x)} \\
&& \leq f(t_{+})g(t_{+})+ f(s_{+})g(s_{+})-f(a_{-})g(a_{-})- f(t_{+})g(s_{+})
 +\sum_{d\in D([a,t])} A(d)\nonumber
\eeq
}

{\bf Proof. } Note first that
\[
 T(f,g; I,J)+T(g,f,I,J) = T(f,g; I,I) +T(f,g; J,J).
 \]
By Theorem \ref{teodi}, we have that
\beqn
T(g,f;I,J) &\geq & T(g,f;J,J) + g_{\min}(I)(\lee{f}(J)-\lee{f}(I)) \\
 & & + \bar{g}(J) (\ree{f}(I) -\ree{f}(J)),
 \eeqn
which  together with the observation that $T(f,g; J,J)= T(g,f; J,J)$ ends the proof.

\eitem

\qitem{Remark}{The prescribed order relation $J\subseteq I$, or, in our particular case, $a\leq s\leq t$, in theorems
\ref{teodi} and \ref{teoco} is important because the L-S integral is sensitive to the order of limits of integration. However, analogous
results can be stated for any other arrangement of $a,s,t$.
}

\eitem

\qitem{Remark}{ It is difficult to state universal necessary conditions for equality in theorems \ref{teodi} and \ref{teoco}.
As an illustration, consider a function $g(x)=\sgn (x-c)$, for a fixed $c\in \R$, and let $f$ be a continuous non-decreasing function.
Then, for $a<c<s<t$, we have that
\[S(f,g; s,t)= \int_{[a,t]} f(x)\di g(x) +\int_{[a,s]} g(x)\di f(x) = f(s) + f(a), \]
which gives equality with both lower and the upper bound. If $a<s<c<t$, then $S(f,g; s,t) = 2f(c)+f(a)-f(s)$, the lower bound is
 $f(s)+f(a)$ and the upper bound is $ 2f(t)+f(a)-f(s)$ and equality with the lower bound occurs iff $f(c)=f(s)$, while the equality with the
 upper bound occurs iff $f(s)=f(t)$.

}

\sect{Examples}

\qitem{Young's inequalities with respect to a measure}{Let $\mu$ be a positive finite
countably additive measure  on Borel subsets of $\R$.
Along the lines of Section \ref{lsip}, we define a right continuous function $G(x)= \mu \{(-\infty, x]\}$. For a non-decreasing function  $f$
we have that
\[  T(f,G ; I,J)=\int_{I} f(x)\di \mu(x) + \int_{J}  G(x)\di f(x) \]

Then we may  observe  Young's inequalities as in  Section \ref{ylowup}. In order to simplify notations, let us assume that $\mu$ is the
probability distribution of a random variable $Y$ on some abstract probability
space $(\Omega, \cF, P)$; then $G$ is the cumulative distribution function
of $Y$, and taking $I=(-\infty,t]$, $J=(-\infty,s]$, we have that
\beq
\label{yim12}
\lefteqn{f(s_{+}) G(t)   +  \sum_{d\in D((-\infty,s])} (f(d)-f(d_{-}))P(Y=d)}\nonumber\\
\lefteqn{  \rule{5em}{0em} \leq \E \left( f(Y) \cdot I_{\{Y\leq t\}}\right) + \int_{[f(-\infty),f(s_{+})]}  G(f^{-1}(u)) \di u}\\
& \leq & f(t_{+})(G(t)-G(s)) +f(s_{+}) G(s) +  \sum_{d\in D((-\infty,t])} (f(d)-f(d_{-}))P(Y=d).\nonumber
 \eeq

This holds for a general probability measure $\mu$. In a purely discrete case  $\mu$ is concentrated  on a countable set of
points $\{y_k\}$, where
we may think of $\{y_k\}$ as being an increasing sequence in $k\in (-\infty, +\infty)$ or $k\in (1,+\infty)$, with
corresponding
probabilities $p_k$.  In this case the Young's functional with $I=(-\infty, y_n]$ and $J=(-\infty, y_m]$, where $m\leq n$, is expressed
by sums as follows.
\[ T(f,G; I,J)= \sum^n f(y_k)p_k +\sum^{m-1} G(y_k) (f(y_{k+1-})-f(y_{k-})) + G(y_m) (f(y_{m+})-f(y_{m-}).\]

Now suppose that we have another probability measure that corresponds to a random variable $X$, with a probability distribution function $F$.
Rewriting (\ref{yim12}) in this setup, $f$ being replaced with $F$, we get

\beqn
\lefteqn{F(s) G(t)   +  \sum_{d\in D((-\infty,s])} P(X=d)P(Y=d)}\nonumber\\
\lefteqn{  \rule{5em}{0em} \leq \E \left( F(Y) \cdot I_{\{Y\leq t\}}\right) + \E \left( G(X) \cdot I_{\{X\leq s\}}\right)}\\
& \leq & F(t)(G(t)-G(s)) +F(s) G(s) +  \sum_{d\in D((-\infty,t])} P(X=d)P(Y=d).
 \eeqn

The case $X=Y$, $F=G$, $s=t=+\infty$ gives a simple formula
\be
\label{efx}
  \E F(X) = \frac{1}{2} +\frac{1}{2}\sum_{d\in D} P^2(X=d),
 \ee
where $D$ is the set of atoms of probability distribution for $X$. Proceeding further, let $F^{-1}$ be any version of inverse; then (\ref{efx})
implies that
\[ F^{-1} (\E F(X)) \geq F^{-1}\left(\frac{1}{2}\right)\in \Med F, \]
where $\Med F$ denotes the median set (closed interval or a singleton), which is defined as the set of all possible values of $F^{-1}(1/2)$.

 By (\ref{geninv}), $F^{-1}\leq \bar{F}^{-1}$, where
\[\bar{F}^{-1} (y) = \inf \{t\; |\; F(t) > y\} \]
 and so, if $m$ is any median of $X$, then
 \[ m \leq \bar{F}^{-1} (\E F(X)) .\]
 This is a new general relation between median and expectation, that seems not to have been recorded in the literature (see for example
 \cite{ghmuk06} or \cite{jenej} for some results and facts about one dimensional medians).

}

\eitem

\qitem{A summation formula}{In this part, the letters $i,j,k,m,n$ will be reserved for integers only.
The following result was proved in \cite{gupta76}: If $f$ is a strictly increasing
positive function on  $[0,+\infty)$, $0<f(1)\leq 1$, then for any positive integer $n$  it holds that
\be
\label{gupsum}
 \sum_{j=1}^n [f(j)] + \sum_{k=1}^{[f(n)]} [f^{-1}(k)] = n[f(n)] + K(1,n),
 \ee
where $K$ is the number of integer points $j\in [1,n]$ such that $f(j)$ is an integer. A generalized version of this
result can be derived using the Lebesgue-Stieltjes integral as follows.

We start with   (\ref{mgenci}) with $g(x)=[x]$ and  replacing $f$ with $[f]$   (where $[x]$ is
the greatest integer  $\leq x$). To be consistent with usual notations, we may use $[f](x)$ for $[f(x)]$; note that
$[f](x_{\pm})$ and $[f(x_{\pm})]$ need not be equal. For simplicity take $a=m$ and $b=n$ where $m$ and $n$ are integers. Then we
have,

\beqn
 \lefteqn{\int_{[m,n]} [f](x)\di [x] +\int_{[m,n]} [x] \di [f](x) =}\\
 & &  n\cdot [f](n_{+}) - (m-1)\cdot [f](m_{-}) + K(m,n),
 \eeqn
where by (\ref{deforc}),
\[
 K(m,n)= \sum_{j=m}^n \big( [f](j)- [f](j_{-})\big) .
 \]

 Now let us fix the version of inverse to be the smallest one:
\[ f^{-1}(y)= \sup\{ t\; |\; f(t)<y\} \]
and also assume that $f$ is right continuous. With these assumptions, we have that
$J(x):=[f](x_{+})-[f](x_{-})= [f](x)-[f](x_{-})\neq 0$ if and only if $f(x)$ is an integer,
and then $J(x)$ equals  the number of integers $k\in ([f](x_{-}),[f](x_{+})] $, or, in another way, the number
of integers $k$ with the property that $f^{-1}(k)=x$. Then
\[\int_{[m,n]} [x] \di [f](x) = \sum_{x\in D([m,n])} [x]([f](x_{+})-[f](x_{-}))= \sum_{k\in ([f](m_{-}),[f](n_{+})]} [f^{-1}(k)] \]
and as an immediate consequence we have that
\beq
\label{sumfo}
\lefteqn{\sum_{j= m}^{n} [f(j)] +\sum_{k=[f](m_{-})+1}^{[f](n_{+})} [f^{-1}(k)]}\\
& =& n \cdot [f](n_{+}) - (m-1)\cdot [f](m_{-}) + K(m,n),\nonumber
\eeq
where $K(m,n)$ is the number of integers $k \in ([f](m_{-}), [f](n_{+})]$ with the property that $f^{-1}(k)$ is also an integer.
This is a generalization of (\ref{gupsum}) for right continuous non-decreasing functions $f$.
The formula could be useful in number theory, for example, to find a numerical value of  $K(m,n)$ for some functions $f$.

}

\eitem

{\bf Acknowledgement. }  The first author acknowledges his partial affiliation to Ra\v cunarski
fakultet, Beograd, Serbia. The work is partially supported by Ministry of Science and Education of Republic of Serbia, projects
III 44006 and 174024.

\medskip

\noindent Department of Applied Mathematics, Faculty of Electrical Engineering, \\
University of Belgrade, P.O. Box 35-54, 11120 Belgrade, Serbia\\
emerkle@etf.rs

\medskip


\noindent Colegiul National "Iancu de Hunedoara", Hunedoara,\\
Str.l Libertatii No.2, Bl. 9, Ap.14,  331032 Romania\\
marinescuds@yahoo.com

\medskip


\noindent Instituto de Mathematica, Universidade Federal do Rio de Janeiro,
Av. Athos da Silveira Ramos 149, C. T. - Bloco C,  21941-909, Rio de Janeiro, Brazil\\
monica@ufrj.br

\medskip


\noindent Colegiul National "Decebal" Deva,\\
Str. Decebal Bl. 8 Ap.10, 330021, Romania\\
mihaimonea@yahoo.com

\medskip


\noindent Colegiul Economic "Emanoil Gojdu", Hunedoara,\\
Str. Viorele No. 4, Bl 10, Ap.10, 331093, Romania \\
stroe@yahoo.com

\end{document}